\documentclass[11pt]{article}    % Enable this line and disable the
\usepackage{graphicx}          % Include this line if your
\usepackage[latin1]{inputenc}
\usepackage[T1]{fontenc}
\usepackage[francais,english]{babel}
\usepackage{amsmath}
\usepackage{amssymb,a4wide}
\newtheorem{theom}{Theorem}
\newtheorem{lm}{Lemma}
\newtheorem{prp}{Proposition}
\newtheorem{rmq}{Remark}
\newcommand{\wmin}{\omega_{*}}
\newcommand{\wmax}{\omega^{*}}
\newcommand{\RR}{{\mathbb R}}
\newcommand{\NN}{{\mathbb N}}
\newcommand{\SSS}{{\mathbb  S}}

\newcommand{\EE}[1]{E\left(#1\right)}
\newcommand{\refeq}[1]{(\ref{#1})}
\newcommand{\dotex}{{\frac{d}{dt}}}

\newcommand{\dvv}[2]{{\frac{\partial^2 #1}{\partial {#2}^2}}}

\begin{document}

\title{Stabilization for an ensemble of  half-spin systems\thanks{KB and PR were partially supported by the ``Agence Nationale de la Recherche'' (ANR),
Projet Blanc C-QUID number BLAN-3-139579. PSPS was partially supported by CNPq --
Conselho Nacional de Desenvolvimento Cientifico e Tecnologico -- Brazil, under grant 308465/2006-7.}} 

\author{ Karine Beauchard\thanks{CMLA, ENS Cachan, CNRS, UniverSud, 61, avenue du Pr\'{e}sident Wilson, F-94230 Cachan, FRANCE. Email: karine.beauchard@cmla.ens-cachan.fr} 
       \and 
       Paulo S\'ergio Pereira da Silva\thanks{Escola Polit\'ecnica da USP - PTC
Cep 05508-900 -- S\~ao Paulo -- SP -- BRAZIL.Email: paulo@lac.usp.br}
       \and              
     Pierre Rouchon\thanks{Mines ParisTech, Centre Automatique et Syst\`{e}mes,
Unit\'{e} Math\'{e}matiques et Syst\`{e}mes,
 60 Bd Saint-Michel, 75272 Paris cedex 06, FRANCE. Email: pierre.rouchon@mines-paristech.fr} 
      }

\maketitle

\begin{abstract}                          % Abstract of not more than 200 words.
Feedback stabilization of  an ensemble of non interacting half spins
described by Bloch equations  is considered.  This system may be seen
as a prototype for infinite dimensional
systems with continuous spectrum. We propose an explicit feedback law that stabilizes
asymptotically the system around a uniform state of spin +1/2 or -1/2.
The proof of the convergence is done locally around the equilibrium in
the $H^1$ topology. This local convergence is shown to be a weak asymptotic
convergence for the $H^1$ topology and  thus a strong convergence for the $C^0$ topology. The proof relies on an adaptation of the LaSalle invariance principle
to infinite dimensional systems.  Numerical simulations illustrate
the efficiency of these feedback laws, even for initial conditions
far from the equilibrium
\end{abstract}


%OOOOOOOOOOOOOOOOOOOOOOOOOOOOOOOOOOOOOOOOOOOOOOOOOOOOOOOOOOOOOOOOOOOOOOOOOOO
\section{Introduction}
%OOOOOOOOOOOOOOOOOOOOOOOOOOOOOOOOOOOOOOOOOOOOOOOOOOOOOOOOOOOOOOOOOOOOOOOOOOO

%---------------------------------------------------------------------------
\subsection{Infinite dimensional systems with continuous spectra}
%---------------------------------------------------------------------------

Most controllability   results available for infinite dimensional bilinear systems are
related to systems with discrete spectra (see for instance,
\cite{beauchard-coron:06} for exact controllability results and
\cite{chambrion-et-al:IHP09,nersesyan:cmp2010} for approximate controllability results).
As far as we know, very few  controllability studies consider
systems admitting a continuous part in their spectra.

In~\cite{mirrahimi-IHP09} an approximate controllability result is given
for a system with mixed  discrete/continuous spectrum:
the Schr\"{o}dinger  partial differential equation of a quantum particle
in an N-dimensional decaying potential  is shown to be approximately controllable
(in infinite time) to the ground bounded state when the  initial state is a
linear superposition of bounded states.

In~\cite{li-khaneja:PRA06,li-khaneja:ieee09} a controllability
notion, called ensemble controllability, is  introduced and  discussed for quantum
systems described by a family of ordinary differential equations (Bloch equations)
depending continuously on a finite number of scalar parameters and with a finite
number of control inputs. Ensemble controllability means that it is possible
to find open-loop controls that compensate for the dispersion in these scalar
parameters: the goal is to simultaneously steer a continuum of systems between
states of interest with the same control input.
Such continuous family of  ordinary differential systems sharing the same
control inputs  can be seen as the prototype of infinite dimensional
systems with purely continuous spectra.

The article~\cite{li-khaneja:ieee09} highlights
the role of Lie algebras and non-commutativity in the design
of a compensating control sequence and  consequently  in the characterization of
ensemble controllability. In \cite{beauchard-et-al:JMP09}, this analysis is completed
by  functional analysis methods developed for infinite
dimensional systems governed by partial differential equations
(see, e.g., \cite{coron:book} for samples of these methods).
Several mathematical answers are given,
with discrimination between approximate and exact controllability,
and finite time and infinite time controllability,
for the Bloch equation. In particular, it is proved that a priori bounded $L^2$-controls
are not sufficient to achieve exact controllability, but unbounded controls
(containing, for example sums of Dirac masses) allow to recover controllability.
For example, it is proved in \cite{beauchard-et-al:JMP09} that the Bloch
equation is approximately controllable to $e_3$, in $H^1$, in finite time,
with unbounded controls. The authors also propose explicit open loop (unbounded) controls
for the local exact controllability to $e_3$ in infinite time.

The goal of this article is to investigate feedback stabilization
of such specific infinite dimensional systems with continuous spectra.
As in~\cite{mirrahimi-IHP09},  the feedback design is based on a Lyapunov
function  closely related to  the  norm  of the state space, a Banach  space.

%---------------------------------------------------------------------------
\subsection{The studied model}
%---------------------------------------------------------------------------

We consider here  an ensemble of non interacting half-spins
in a static field $(0,0,B_0)^T$ in $\mathbb{R}^3$,
subject to a transverse radio frequency field
$(\tilde{u}(t),\tilde{v}(t),0)^T$ in $\mathbb{R}^3$
(the control input). The ensemble of half-spins is
described  by the magnetization vector $M\in\mathbb{R}^3$ depending on time $t$
but also on the Larmor frequency $\omega=-\gamma B_0$ ($\gamma$ is the gyromagnetic
ratio). It obeys to  the Bloch equation:
\begin{equation}
 \label{dyn:eq1}
  \frac{\partial M}{\partial t}(t, \omega) =
(\tilde u(t) e_1 + \tilde v(t) e_2 + \omega e_3) \times
  M(t, \omega),
\end{equation}
where
$-\infty < \wmin < \wmax < +\infty$,
$\omega\in(\wmin,\wmax)$,
$(e_1,e_2,e_3)$ is the canonical basis of $\mathbb{R}^3$,
$\times$ denotes the vector  product on $\mathbb{R}^3$.
The equation (\ref{dyn:eq1}) is an infinite dimensional bilinear control system in which
\begin{itemize}
\item the state is the function $M$, with, for every $\omega \in (\wmin,\wmax)$,
$M(t,\omega)\in  \mathbb{S}^{2}$, the unit sphere of $\mathbb{R}^3$,
\item the two control inputs  $\tilde{u}$ and $\tilde{v}$ are real valued.
\end{itemize}
It must be stressed that $\tilde u(t)$ and
$\tilde v(t)$ are common controls for all the members of the
ensemble, and they cannot depend on $\omega$.
In coordinates $M=(x,y,z)$, the Bloch equation may be written
 \begin{equation}
  \label{dyn:eq}
\left\lbrace \begin{array}{rcl}
 \dot x & = & - \omega y + \tilde v z,\\
 \dot y & = &\omega x - \tilde u z, \\
 \dot z & = & -\tilde v x + \tilde u y,
\end{array} \right.
\end{equation}
where $\dot x$ stands for $\partial x/\partial t$.

Now, let us precise the definition of a solution associated to Dirac controls.
When $\tilde{u}, \tilde{v} \in L^1_{loc}(\mathbb{R})$, then,
for every initial condition $M_0 \in L^2((\omega_*,\omega^*),\mathbb{R}^3)$,
the equation (\ref{dyn:eq1}) has a weak solution in the usual sense:
$M \in C^0([0,+\infty),L^2((\wmin,\wmax),\mathbb{R}^3))$.
Let us give the definition of a solution associated to a Dirac control. In the sequel $\delta_a(t)$ stands for the Dirac function $\delta(t-a)$.
If $\tilde{u}=\alpha \delta_a (t)+u_{\sharp}$ and $\tilde{v}=v_{\sharp}$
where $u_{\sharp}, v_{\sharp} \in L^1_{loc}(\mathbb{R})$,
$\alpha>0$ and $a \in (0,+\infty)$, then
the solution is the classical solution on $[0,a)$ and $(a,+\infty)$,
it is discontinuous at the point $a$, with the explicit discontinuity
$$M(a^+,\omega)=\left( \begin{array}{ccc}
1 & 0            &  0           \\
0 & \cos(\alpha) & -\sin(\alpha)\\
0 & \sin(\alpha) &  \cos(\alpha)
\end{array}\right) M(a^-,\omega).$$
This definition corresponds to the limit, as $\epsilon \rightarrow 0^+$
of solutions associated to $\tilde{u}(t)=(\alpha/\epsilon) 1_{[a,a+\epsilon]}(t)$.

Formally, the spectrum of the operator $A$ defined by
$$(AM)(\omega):=\omega e_3 \times M(\omega)$$
is $-i (\wmin,\wmax) \cup i (\wmin,\wmax)$: for every $\omega_\sharp \in (\wmin,\wmax)$,
the eigenvector associated to $\pm i \omega_\sharp$ is
$(1, \mp i, 0)^T \delta_{\omega_\sharp}(\omega)$.
Thus the Bloch equation (\ref{dyn:eq1}) is a prototype of infinite dimensional system
with continuous spectrum.

%---------------------------------------------------------------------------
\subsection{Outline}
%---------------------------------------------------------------------------

The goal of this article it to propose a first answer to the following question.
\vspace{\baselineskip}

{\bf Local Stabilization Problem.}
Define an explicit control law
$(\tilde{u},\tilde{v})=(\tilde{u}(t,M),\tilde{v}(t,M))$
and a neighborhood $U$ of $-e_3$
(in some space of functions $(\wmin,\wmax) \rightarrow \SSS^2$
to be determined) such that, given any initial condition
$M^0 \in U$, the solution of the closed loop system
is defined for every $t \in [0,+\infty)$,
is unique and converges to $-e_3$, when $t \rightarrow + \infty$,
uniformly with respect to $\omega \in (\wmin,\wmax)$.
\vspace{\baselineskip}

Section~\ref{H1:sec}  is devoted to control design and closed-loop simulations:
the  feedback law is the sum of a Dirac comb and
a time-periodic feedback law based on a Lyapunov function;
Proposition~\ref{Solutions} proved in Appendix A guarantees that
the closed-loop initial value problem is always well defined;
simulations illustrate the convergence rates observed for
an initial state formed by a quarter of the equator  on the Bloch sphere.
In section~\ref{main:sec} we state and prove the main convergence result,
Theorem~\ref{MainThm}: the closed-loop convergence towards the constant profile $M(\omega)=-e_3$
is shown to be  local  and weak  for the  $H^1$ topology on $M$.
The obstruction to global stabilization is also discussed: it is based on an explicit description of the Lasalle invariant set.
Some concluding remarks are gathered in section~\ref{conclusion:sec}.

%OOOOOOOOOOOOOOOOOOOOOOOOOOOOOOOOOOOOOOOOOOOOOOOOOOOOOOOOOOOOOOOOOOOOOOOOO
\section{Lyapunov $H^1$ approach}{\label{H1:sec}
%OOOOOOOOOOOOOOOOOOOOOOOOOOOOOOOOOOOOOOOOOOOOOOOOOOOOOOOOOOOOOOOOOOOOOOOOO

%--------------------------------------------------------------------------
\subsection{The impulse-train control}
%--------------------------------------------------------------------------

It is proved in~\cite{beauchard-et-al:JMP09} that
controls containing sums of Dirac masses are crucial to achieve the controllability
of the Bloch equation. In view of the controls used in this reference,
it is natural to consider a control with the following ``impulse-train'' structure
\begin{equation}\label{uv:eq}
\tilde {u} = u + \sum_{k=1}^{+\infty} \pi~\delta_{k\!T}(t), \quad
\tilde{v} = (-1)^{\EE{\frac{t}{T}}} v
\end{equation}
for some period $T>0$, which is fixed in all the article
($\EE{\gamma}$ denotes the integer part of the real number $\gamma$).
The new controls $u$ and $v$ belong to
$L^1_{loc}(\mathbb{R})$. Then, after each impulse that is
applied at time $t=kT$, $x$ remains unchanged, but  $y$
and $z$ are moved to their opposites, that is
$$(x,y,z)(kT^+)=(x,-y,-z)(kT^-).$$
The resulting state diffeomorphism
 \begin{equation}
  \label{e:identification}
(x,y,z) \mapsto (\mathcal{X}=x,\mathcal{Y}=-y,\mathcal{Z}=-z)
 \end{equation}
 transforms~\eqref{dyn:eq} into
\begin{eqnarray*}
 \dot{ \mathcal{X}}(t) & = & \omega \mathcal{Y}(t) - \tilde v(t) \mathcal{Z}(t),\\
 \dot {\mathcal{Y}}(t) & = & -\omega \mathcal{X}(t) - \tilde u \mathcal{Z}(t), \\
 \dot {\mathcal{Z}}(t) &= & \tilde v(t) \mathcal{X}(t) + \tilde u(t) \mathcal{Y}(t).
\end{eqnarray*}
Let $\varsigma=(-1)^{\EE{\frac{t}{T}}}$. Considering the following change of variables
$$(x,y,z)(t,\omega)
\mapsto
\Big( x(t,\omega) , \zeta(t)y(t,\omega) , \zeta(t)z(t,\omega) \Big)$$
one gets the following dynamics
\begin{equation}  \label{dynS:eq}
\left\lbrace \begin{array}{l}
 \dot{x}  =  - \varsigma \omega y + v z,\\
 \dot{y}  = \varsigma  \omega x -  u z,\\
 \dot{z}  =  - v x+  u y,
\end{array} \right.
\end{equation}
with the new control $( u,  v)$ as in (\ref{uv:eq}).
It is as if, between $[kT,(k+1)T]$
and $[(k+1)T,(k+2)T]$, one is changing the sign of $\omega$, but the
solution, after the identification \refeq{e:identification}, remains
continuous in $t$ (but not differentiable in $t$ at the instants $t = k T,
k \in \NN$). In other words, the application of the impulses at $t =
k T$ changes the sense of rotation of the null input solution. One
would expect that  this impulse-train control is reducing the
average dispersion of the solution. Roughly speaking, the dispersion observed for the open-loop system~\eqref{dyn:eq}  with $(\tilde u,\tilde v)$ as  input  is strongly reduced and almost canceled for  the open-loop system~\eqref{dynS:eq} with $(u,v)$ as  input.

%------------------------------------------------------------------------
\subsection{Heuristics of the Lyapunov-like control}
%------------------------------------------------------------------------

Now let $Z(t, \omega)$ and $\Omega(t)$ defined by
\[
Z = x + i y, \quad \Omega = v - i u
\]
where $x, y, z$ refer to the transformed dynamics \refeq{dynS:eq}.
Then one may write \refeq{dynS:eq} in the form
\[
 \left\{
 %%%%
 %%%%
\begin{array}{rcl}
\dot{Z}(t , \omega) & = & i \varsigma(t) \omega Z(t, \omega) +
 \Omega(t) z(t, \omega),\\
 \dot{z}(t, \omega) & = & - \Re \left[ \Omega(t) \overline{Z(t, \omega)}\right],
\end{array}
 \right.
\]
where $\Re(\xi)$ (resp. $\overline{\xi}$) denotes the real part
(resp. the complex conjugate number) of a complex number $\xi \in \mathbb{C}$.
The following transformation
\[
 \tilde {Z}(t,\omega) = Z(t, \omega) e^{- i \omega \int_{0}^t \varsigma(\tau) d\tau}
\]
converts the system into the driftless form
\begin{equation}
\label{e:complex:form}
 \left\{ \begin{array}{rcl}
\dot Z (t,\omega) & = & \Omega(t) z(t,\omega) e^{- i \omega
\int_{0}^t \varsigma},\\
\dot{z}(t, \omega) & = & - \Re \left[ \Omega(t) \overline{Z(t,
\omega)} e^{- i \omega \int_{0}^t \varsigma} \right],
\end{array}
 \right.
\end{equation}
where, for notation simplicity, one lets $Z(t, \omega)$ stand for
$\tilde Z(t,\omega)$, and one lets $\int_{0}^t \varsigma$ stand for
$\int_{0}^t \varsigma(\tau) d\tau$.

For the moment one shall assume that the input $\Omega(t)$ will be
chosen in such a way that the solution $(Z(t,\omega), z(t, \omega))$ of
\refeq{e:complex:form} does exist, it is unique and it is regular
enough in a way that one may consider that the derivatives
$Z^\prime(t,\omega) = \frac{\partial{Z}}{\partial \omega}(t,\omega)$
and $z^\prime(t,\omega) = \frac{\partial{z}}{\partial
\omega}(t,\omega)$ exist almost everywhere and they are solutions
of the differential equation that is obtained by differentiation of
\refeq{e:complex:form} with respect to $\omega$, namely
 \begin{equation} \label{e:prime}
 \left\lbrace \begin{array}{l}
 \dot{Z}^\prime  =  \Omega \left\{ \left[z^\prime  - i \left(\int_{0}^t
 \varsigma\right) z \right] e^{-i \omega \int_{0}^t
 \varsigma}  \right\},\\
 \dot{z}^\prime  =  - \Re \left\{ \Omega \left[ \overline{Z}^\prime
- i  \left(\int_{0}^t
 \varsigma\right) \overline{Z} \right]e^{-i \omega \int_{0}^t
 \varsigma} \right\},
 \end{array} \right.
 \end{equation}
where $\dot{Z}^\prime$ stands for $\frac{\partial}{\partial t}
{Z}^\prime$, and $\dot{z}^\prime$ stands for
$\frac{\partial}{\partial t} {z}^\prime$.

Now consider the following Lyapunov-like functional:
\begin{equation}~\label{lyap:eq}
 \mathcal{L} =  \frac{1}{2} \int\limits_{\wmin}^{\wmax}
\left\{  | Z^\prime |^2 + (z^\prime)^2  + 2 G (z+1)  \right\} d\omega
\end{equation}
where $G$ is a positive real number and $Z(t, \omega),
Z^\prime(t,\omega)$ and $z^\prime(t, \omega)$ refer to the
solutions respectively of (\ref{e:complex:form}) and
(\ref{e:prime}). One may write
\begin{equation} \label{dL/dt=}
 \dotex \mathcal{L}(t) = \Re\left(
\int\limits_{\wmin}^{\wmax}
\left\{  \bar{Z}^\prime \dot{Z}^\prime  + z^\prime \dot{z}^\prime +  G \dot{z} \right\} d\omega\right)
\end{equation}
and so, taking into account (\ref{e:complex:form}) and
(\ref{e:prime}), the fact that
$\Omega(t)$ does not depend on $\omega$, one gets
\begin{equation} \label{dL/dt=bis}
 \dotex \mathcal{L}(t) =
 \Re \left[ \Omega(t) H(t) \right]
\end{equation}
where
$$H(t) :=   \int\limits_{\wmin}^{\wmax} \Big\{
i \left(\int\limits_{0}^t \varsigma\right)
\left( \bar Z z^\prime - \bar{Z}^\prime z  \right)
 -  G \bar Z
\Big\} e^{-i \omega \int_{0}^t \varsigma} d\omega.$$
Hence one may take $\Omega(t) = - K_p \bar H(t)$, where $K_p$ is a positive real number, obtaining
\begin{equation} \label{e:control:law}
\Omega(t) =  K_p \int\limits_{\wmin}^{\wmax} \Big\{
i \left(\int\limits_{0}^t \varsigma \right)  \left( Z z^\prime - {Z}^\prime z \right)
+ G Z   \Big\} e^{i \omega \int_{0}^t \varsigma} d\omega.
\end{equation}
It follows that
$$\frac{d \mathcal{L}}{dt} (t) = - \frac{1}{K_p} | \Omega(t) |^2.$$

Consider the system \refeq{e:complex:form} in closed loop with the control law
\refeq{e:control:law}, thereby called by \emph{closed loop system}.
The state space of this system is $H^1((\wmin, \wmax), \RR^3)$,
which is the set of functions $f \in L^2(\wmin,\wmax)$ such that
the distributional derivative $f'$ belongs to $L^2(\wmin,\wmax)$.
This space, equipped with the norm
$$\|f\|_{H^1}:=\left( \int\limits_{\wmin}^{\wmax}
|f'(\omega)|^2 + |f(\omega)|^2 d\omega \right)^{1/2}$$
is a Banach space. In other words,  closed-loop system~(\ref{e:complex:form})-(\ref{e:control:law}) may be considered
to be a differential equation of the form
$$\left\lbrace \begin{array}{lll}
\dot{M}(t) & = & F(t, M(t))\\
 M(0) & = & M^0 \in H^1((\wmin, \wmax), \SSS^2)
\end{array} \right.$$
where $F(t, M)$ is a continuous map
 \[ F : \RR \times H^1((\wmin, \wmax), \RR^3) \rightarrow
H^1((\wmin, \wmax), \RR^3).
 \]
Moreover, $F$ is locally integrable ($L^1_{loc}$) and periodic in $t$
and locally Lipschitz in $M$.
Using the same ideas as in the proof of the
Cauchy-Lipschitz (Picard-Lindel\"of) theorem,
we get local (in time) solutions in $C^0([0,T],H^1)$.
From the contruction of the feedback law,
finite time blow up in $H^1$ is impossible,
thus solutions are global in time.
Precisely, we have the following result,
whose proof is detailed in Appendix A.
\begin{prp} \label{Solutions}
 For every initial condition $M^0 \in H^1((\wmin, \wmax), \SSS^2)$,
the closed loop system (\ref{e:complex:form}), (\ref{e:control:law})
has a unique solution $M \in C^1_{pw} \left([0, \infty), H^1\left( (\wmin, \wmax), \mathbb{S}^2 \right) \right)$
such that $M(0)=M^0$.
\end{prp}
In this statement, the space $C^1_{pw} \left([0, \infty), H^1\left( (\wmin, \wmax), \SSS^2 \right) \right)$
is made of functions\\ $M \in C^0 \left([0, \infty), H^1\left( (\wmin, \wmax), \SSS^2 \right) \right)$ such that
$M \in C^1((kT,(k+1)T),H^1\left( (\wmin, \wmax), \SSS^2 \right))$, for every $k \in \mathbb{N}$;
their derivative $\partial M / \partial t$ is possibly discontinuous at $t=k T$,
but it has finite limits in $H^1\left( (\wmin, \wmax), \mathbb{S}^2 \right)$
when $t \rightarrow (kT)^+$ and $t \rightarrow (kT)^-$.

%----------------------------------------------------------------------------------
\subsection{Closed-loop simulations}
%----------------------------------------------------------------------------------

We assume here $\wmin=0$, $\wmax=1$ and we solve numerically the
$T$-periodic system~\eqref{dynS:eq} with the $T$-periodic feedback
law~\eqref{e:control:law} where $Z = x + \imath y$ and  $\Omega = v
- \imath u $. The parameters are $T=2\pi/(\wmax-\wmin)$, $K_p=1$,
$G=T^2/20$. The simulation is for $t\in[0, T_f]$, $T_f=20 T$. The
$\omega$-profile $[\wmin,\wmax]\ni \omega \mapsto
(x(t,\omega),y(t,\omega),z(t,\omega))$ is discretized $\{0,\ldots,
N\}\ni k \mapsto (x_k(t),y_k(t),z_k(t))$ with a regular mesh of step
$\epsilon_N = \frac{\wmax-\wmin}{N}$ and $N=100$:
$(x_k(t),y_k(t),z_k(t))$ is  then an approximation of
$(x(t,k\epsilon_N),y(t,k\epsilon_N),z(t,k\epsilon_N))$. We have
checked  that the closed-loop simulations are almost identical for
$N=100$ and $N=200$.  In the feedback law~\eqref{e:control:law}, the
integral versus $\omega$ is computed assuming that $(x,y,z)$ and
$(x^\prime,y^\prime,z^\prime)$ are constant over
   {\small $](k-\frac{1}{2})\epsilon_N,(k+\frac{1}{2})\epsilon_N[$},
   their values being $(x_k,y_k,z_k)$ and {\small $\left(\frac{x_{k+1}-x_{k-1}}{2\epsilon_N},
   \frac{y_{k+1}-y_{k-1}}{2\epsilon_N},\frac{z_{k+1}-z_{k-1}}{2\epsilon_N}\right)$.}
The obtained differential  system is  of  dimension $3(N+1)$. It is
integrated via an explicit Euler scheme with a step size $h= T/800$.
We have tested that $h=T/1600$ yields to almost the same numerical
solution at $t=T_f=20T$. After each time-step the new values of
$(x_k,y_k,z_k)$ are  normalized to remain in $\SSS^2$.

Figures~\ref{fig:LyapControl} and~\ref{fig:BeginEnd}
summarize the main convergence issues when the initial
$\omega$-profiles of $(x,y,z)\in\SSS^2$ are  $z(0,\omega)= 0$,
$x(0,\omega)=\cos\left(\tfrac{\pi}{2}\omega+\tfrac{\pi}{4}\right)$ and
$y(0,\omega)=\sin\left(\tfrac{\pi}{2}\omega+\tfrac{\pi}{4}\right)$. The convergence
speed  is rapid  at the beginning and tends to decrease at the end. We start with ${\mathcal L}(0)\approx 3.24$ and $\sup_{\omega\in[\wmin,\wmax]} \|M(0,
\omega) + e_3\|\approx 1.41$.
We get ${\mathcal L}(20T)\approx 0.38$ and $\sup_{\omega\in[\wmin,\wmax]} \|M(20T,
\omega) + e_3\|\approx 0.27$. The control problem is  quite  hard due to the fact that one has a
continuous spectrum, that is, an infinite ensemble of systems with a
common control input $\Omega(t)$. Hence, as time increases, the
control must fight against the dispersion of the solutions $M(t,
\omega)$ for different values of $\omega$. Simulations (not
presented here) on  much longer times until $1000T$  and with the same
initial conditions and parameters  always yield to smaller final
value for  the Lyapunov function: we get ${\mathcal L}(1000T)\approx 0.02$ and $\sup_{\omega\in[\wmin,\wmax]} \|M(1000T,
\omega) + e_3\|\approx 0.04$.
This is a strong indication of asymptotic converge of the profile
$M(t,\omega)$ toward $-e_3$, even if the convergence speed seems to
be very slow. This numerically observed convergence is confirmed by
Theorem~\ref{MainThm} here below.

\begin{figure}[bht]
  \centerline{\includegraphics[width=\textwidth]{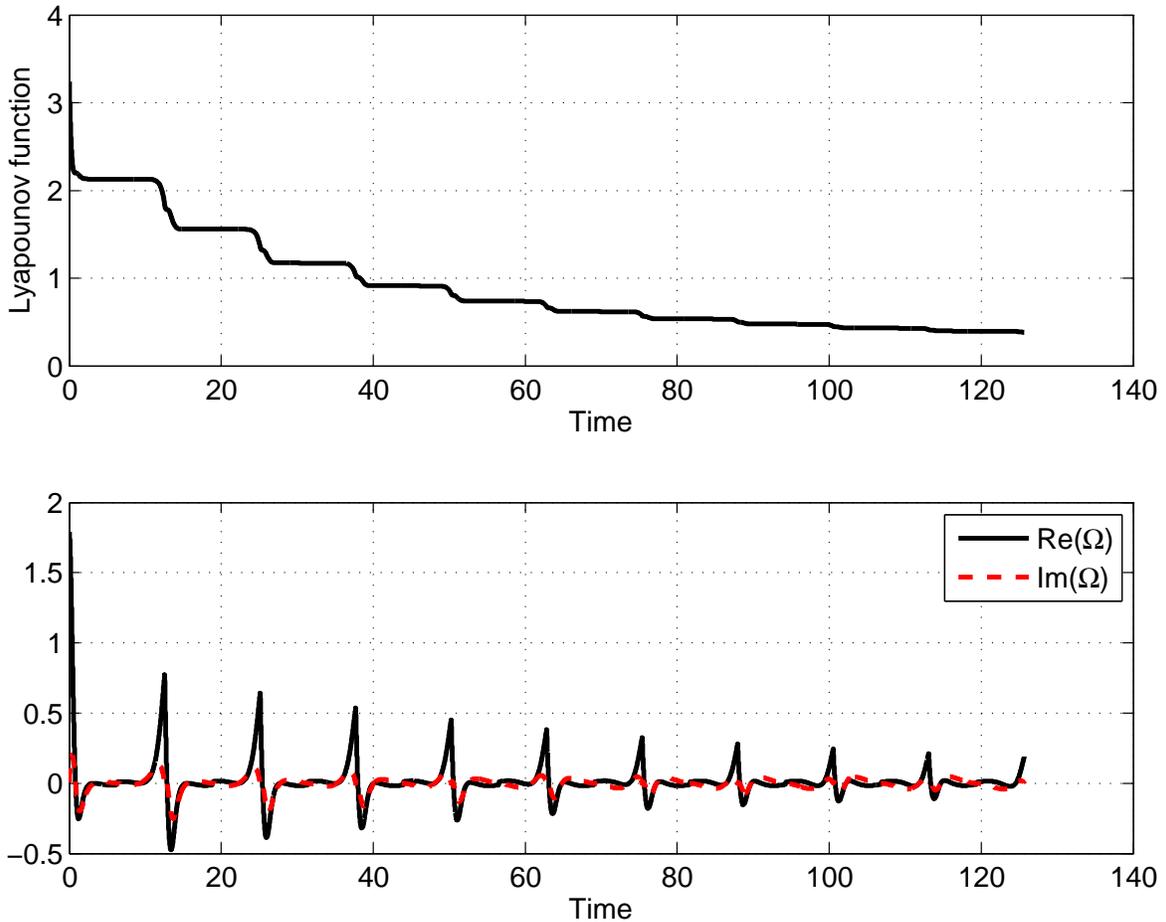}}
  \caption{Lyapunov function $\mathcal{L}(t)$ defined by~\eqref{lyap:eq}
  and  the closed-loop control $\Omega(t)$ defined by~\eqref{e:control:law}
  .}\label{fig:LyapControl}
\end{figure}

\begin{figure}[bht]
  \centerline{\includegraphics[width=\textwidth]{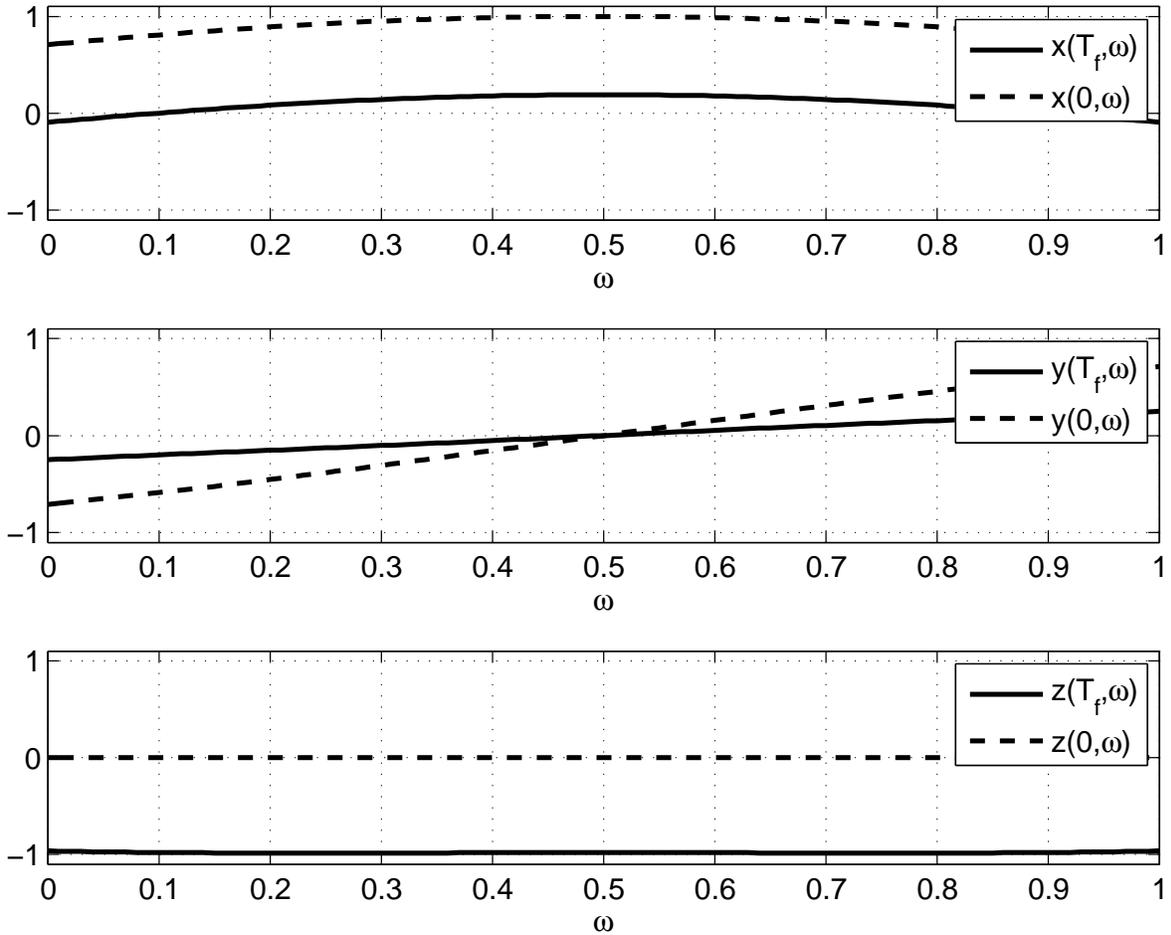}}
  \caption{Initial ($t=0$) and final ($t=T_f$) $\omega$-profiles for
  $x$, $y$ and $z$ solutions of the closed-loop system~\eqref{dynS:eq}
  with the $T$-periodic feedback~\eqref{e:control:law}.}\label{fig:BeginEnd}
\end{figure}

%OOOOOOOOOOOOOOOOOOOOOOOOOOOOOOOOOOOOOOOOOOOOOOOOOOOOOOOOOOOOOOOOOOOOOOO
\section{Main Result} \label{main:sec}
%OOOOOOOOOOOOOOOOOOOOOOOOOOOOOOOOOOOOOOOOOOOOOOOOOOOOOOOOOOOOOOOOOOOOOOO

%----------------------------------------------------------------------
\subsection{Local stabilization}
%---------------------------------------------------------------------

The main result of this paper shows that the control law
\refeq{e:control:law} is a solution of the local stabilization problem stated at the end of the introduction.
\begin{theom} \label{MainThm}
Consider system~\eqref{e:complex:form} with the feedback law~\eqref{e:control:law}.
There exists $\delta' >0$ such that, for every
$M_0 \in H^1((\wmin,\wmax), \SSS^2)$ with $\|M_0 + e_3\|_{H^1} \leq \delta'$,
$M(t,\omega)$ converges  weakly in $H^{1}$ to  $-e_3$
when $t \rightarrow + \infty$.
In particular, as the injection of $H^1$ in $C^0$ is compact,
$M(t,\omega)$ converges to $-e_3$ when $t \rightarrow + \infty$
uniformly with respect to $\omega \in (\wmin,\wmax)$
(convergence in the sup norm of $C^0$).
\end{theom}

The proof of Theorem~\ref{MainThm} relies on an adaptation
of the LaSalle invariance principle to infinite dimensional systems.
This principle is a powerful tool
to prove the asymptotic stability of an equilibrium of
a dynamical system in finite dimension: one just needs to
check that the invariant set coincides with the target \cite{MMPRGT}.
For infinite dimensional systems, the use of the LaSalle invariance
principle is more delicate (because closed and bounded subsets are not necessarily compact).
Roughly speaking, there are 2 adaptation strategies:
\begin{itemize}
\item either one accepts a weaker result:
approximate stabilization \cite{KB-MM,mirrahimi-IHP09}
or weak stabilization \cite{BS};
in this case, one may only need to ensure that the invariant set is reduced to the target,
\item or one wants strong stabilization;
in this case, one needs an additional compactness property for the trajectories of
the closed loop system \cite{JMC-BAN},
or a strict Lyapunov function \cite{JMC-BAN-GB}.
\end{itemize}
This article concerns the first strategy: we prove weak stabilization.

The first step of our proof consists in checking that, locally,
the invariant set is $\{-e_3\}$.

\begin{prp} \label{Prop:Inv}
There exists $\delta>0$ such that,
for every $M^0 \in H^1((\wmin,\wmax),\SSS^2)$ with
$\|M^0+e_3\|_{H^1} < \delta$,  the map
$t \mapsto \mathcal{L}(t)$ is constant on $[0,+\infty)$
if and only if $M^0=-e_3$.
\end{prp}

\noindent \textbf{Proof:} Let us assume that $\mathcal{L}(t)$ is constant.
Then, $\Omega(t) = 0$, $Z(t,\omega)=Z_0(\omega)$ and $z(t,\omega)=z_0(\omega)$.
We deduce from (\ref{e:control:law}) that
$$\int\limits_{\wmin}^{\wmax} \Big\{
 i  t ( Z_0' z_0 -  Z_0 z_0') - G Z_0
\Big\} e^{i \omega t} d\omega = 0, \forall t \in [0,T].$$
Considering the power series expansion versus $t$ of the left hand side, we get
\begin{equation} \label{CNS-inv}
\int\limits_{\wmin}^{\wmax} \left[
P' \Big( Z_0' z_0 - Z_0 z_0' \Big)
- G P Z_0 \right] d\omega = 0, \forall P \in \mathbb{C}[\omega].
\end{equation}
Polynomials are dense in $H^1(\omega_*,\omega^*)$, thus,
the previous inequality holds for every $P \in H^1(\wmin,\wmax)$.
In particular, with $P(\omega)=\overline{Z_0(\omega)}$, we get
$$
\int_{\omega_*}^{\omega^*} \left[
 \overline{Z_0'} \Big( Z_0' z_0 - Z_0 z_0' \Big)
- G \overline{Z_0} Z_0
\right] d\omega = 0.
$$
Summing this equality with the following left hand side, we deduce that
$$\begin{array}{ll}
& \int_{\omega_*}^{\omega^*}
\Big[ |Z_0'|^2 + G |Z_0|^2 \Big] d\omega
\\ = &
\int_{\omega_*}^{\omega^*} \left[
(1+z_0) |Z_0'|^2 - \overline{Z_0'} Z_0 z_0'
\right] d\omega
\\  \leqslant &
\|1+z_0\|_{L^\infty} \|Z_0'\|_{L^2}^2
+ \|Z_0'\|_{L^2} \| Z_0\|_{L^\infty} \|z_0'\|_{L^2}.
\end{array}$$
Thanks to the continuity of the embedding
$H^1(\wmin,\wmax) \subset L^\infty(\wmin,\wmax)$,
there exists a constant $\mathcal{C}>0$ (independent of $M$) such that
$$\| Z_0 \|_{L^\infty}
\leqslant \mathcal{C} \left(
\int_{\wmin}^{\wmax} ( |Z_0'|^2 + G |Z_0|^2 ) d\omega
\right)^{1/2}.$$
Therefore, we get
$$\begin{array}{ll}
& \int_{\omega_*}^{\omega^*}
\Big[ |Z_0'|^2 + G |Z_0|^2 \Big] d\omega
\\ \leqslant &
\Big( \|1+z_0\|_{L^\infty} + \mathcal{C} \|z_0'\|_{L^2} \Big)
\int_{\omega_*}^{\omega^*}
\Big[ |Z_0'|^2 + G |Z_0|^2 \Big] d\omega.
\end{array}$$
There exists $\delta>0$ such that, for every
$M_0 \in H^1((\wmin,\wmax),\SSS^2)$ with
$\| M_0+e_3 \|_{H^1} < \delta$, we have
$$\|1+z_0\|_{L^\infty} + \mathcal{C} \|z_0'\|_{L^2} < 1.$$
If $\mathcal{L}$ is constant along the trajectory
associated to such an initial condition $M_0$,
then, the previous argument shows that $Z_0=0$,
thus $M_0=-e_3$. $\hfill \Box$
\\
\\

For the proof of Theorem \ref{MainThm},
we need the continuity with respect to initial conditions,
of the solutions of the closed loop system (\ref{e:complex:form}), (\ref{e:control:law}),
for the $H^{\frac{1}{2}}(\wmin,\wmax)$-topology. This space is defined by interpolation
between $L^2(\wmin,\wmax)$ and $H^1(\wmin,\wmax)$ and we have a compact injection
$H^1(\wmin,\wmax) \rightarrow H^{\frac{1}{2}}(\wmin,\wmax)$.
We also use the space $H^{-\frac{1}{2}}(\wmin,\wmax)$,
which is the dual space of $H^{\frac{1}{2}}(\wmin,\wmax)$ for the $L^2$-scalar product.
First, let us recall the following Lemma.

\begin{lm} \label{Lem:c1}
There exists $c_1>0$ such that,
for every $\varphi \in H^{\frac{1}{2}}(\wmin,\wmax)$
and for every $\alpha \in [0,T]$,
the map $\omega \mapsto \varphi(\omega) e^{i \alpha \omega}$
belongs to $H^{\frac{1}{2}}(\wmin,\wmax)$ and satisfies
$$\| \varphi(\omega) e^{i \alpha \omega} \|_{H^{\frac{1}{2}}}
\le c_1 \| \varphi \|_{H^{\frac{1}{2}}}.$$
\end{lm}

\noindent \textbf{Proof:} We have
$$\| \varphi(\omega) e^{i \alpha \omega} \|_{L^2} = \| \varphi \|_{L^2},
\forall \varphi \in L^2(\wmin,\wmax),$$
and, for every $\varphi \in H^1(\wmin,\wmax)$,
$$\begin{array}{ll}
\| \varphi(\omega) e^{i \alpha \omega} \|_{H^1}^2
& =
\int\limits_{\wmin}^{\wmax}
| \varphi'(\omega) + i \alpha \varphi(\omega)|^2
+ |\varphi(\omega)|^2 d\omega
\\ & \leqslant \int\limits_{\wmin}^{\wmax}
2 | \varphi'(\omega)|^2 + (2\alpha^2+1) |\varphi(\omega)|^2 d\omega
\end{array}$$
thus we get the conclusion with, for example, $c_1:=(2T^2+2)^{1/4}$
by interpolation. $\hfill \Box$

\begin{prp}
There exists $\delta'>0$ such that,
for every $(M_n^0)_{n \in \mathbb{N}} \in H^1((\wmin,\wmax),\SSS^2)^{\mathbb{N}}$,
$M_\infty^0 \in H^1((\wmin,\wmax),\SSS^2)$ satisfying
\begin{itemize}
\item $\|M_n^0+e_3\|_{H^1} < \delta', \forall n \in \mathbb{N}$,
\item $M_n^0 \rightharpoonup M_\infty^0$ weakly in $H^1$ when $n \rightarrow + \infty$,
\item $M_n^0 \rightarrow M_\infty^0$ strongly in $H^{\frac{1}{2}}$ when
$n \rightarrow + \infty$,
\end{itemize}
the solutions $M_n(t,\omega)$, $M_\infty(t,\omega)$ of the closed loop system
associated to these initial conditions satisfy the following convergences,
when $n \rightarrow + \infty$, for every $t \in [0,+\infty)$,
$$M_n(t) \rightarrow M_\infty(t)
\text{ strongly in } H^{\frac{1}{2}}
\text{ and }
\Omega_n(t) \rightarrow \Omega_\infty(t).$$
\end{prp}

\noindent \textbf{Proof:} First, let us emphasize that $\sqrt{\mathcal{L}(M)}$ and
$\|M+e_3\|_{H^1}$ are equivalent norms on a small enough
$H^1((\wmin,\wmax),\SSS^2)$-neighborhood of $-e_3$: there exists
$\eta, c_*,c^* > 0$ such that, for every
$M \in H^1((\wmin,\wmax),\SSS^2)$ with $\|M+e_3\|_{H^1}<\eta$, we have
\begin{equation} \label{norm-equiv}
c_* \sqrt{\mathcal{L}(M)}
\le \| M+e_3 \|_{H^1}
\le c^* \sqrt{\mathcal{L}(M)}.
\end{equation}
Now, let $\delta':=\min\{\delta c_*/c^*,\eta\}$,
where $\delta$ is as in Proposition \ref{Prop:Inv}.
Thanks to the monotonicity of $\mathcal{L}$, we have,
for every $t \in [0,+\infty)$,
$$\begin{array}{ll}
\|M_n(t)+e_3\|_{H^1}
& \le c^* \sqrt{\mathcal{L}(M_n(t))}
\\ & \le c^* \sqrt{\mathcal{L}(M_n^0)}
\\ & \le \frac{c^*}{c_*} \|M_n^0+e_3\|_{H^1} < \frac{c^* \delta'}{c_*} \le \delta.
\end{array}$$

We have
$$\begin{array}{l}
\| M_n(t)-M_\infty(t)\|_{H^{\frac{1}{2}}} \le
\| M_n^0 - M_\infty^0\|_{H^{\frac{1}{2}}}
\\ + \int_0^t \| F(s,M_n(s)) - F(s,M_\infty(s)) \|_{H^{\frac{1}{2}}} ds.
\end{array}$$
Let us prove the existence of $C>0$ such that,
for every $M,\tilde{M} \in H^1(\wmin,\wmax)$ satisfying $\|M+e_3\|_{H^1} < \delta$,
we have
$$\| F(s,M) - F(s,\tilde{M}) \|_{H^{\frac{1}{2}}}
\le C \|M-\tilde{M}\|_{H^{\frac{1}{2}}},
\forall s \in \mathbb{R}.$$
Then, we will conclude the proof thanks to the Gronwall Lemma.
Let us work, for example, on the third component of $F$:
$$F_3(t,M)=-Re\left[
\Omega(t) \overline{Z} e^{-i \omega \int_0^t \zeta}
\right],$$
where $\Omega$ is defined by (\ref{e:control:law}).
We have
$$\begin{array}{ll}
& \| F_3(t,M)-F_3(t,\tilde{M}) \|_{H^{\frac{1}{2}}}
\\ \le &
| \Omega(t) - \tilde{\Omega}(t) |
\| \overline{Z} e^{-i \omega \int_0^t \zeta} \|_{H^{\frac{1}{2}}}
\\ & + | \tilde{\Omega}(t) |
\| \overline{(Z-\tilde{Z})} e^{-i \omega \int_0^t \zeta} \|_{H^{\frac{1}{2}}}
\\  \le &
| \Omega(t) - \tilde{\Omega}(t) |  c_1 \| Z \|_{H^{\frac{1}{2}}}
+ K c_1 \| Z-\tilde{Z} \|_{H^{\frac{1}{2}}}
\end{array}$$
where $c_1$ is as in the previous Lemma and $K=K(\delta)$.
It is sufficient to prove the existence of a constant $C>0$ such that,
for every $M,\tilde{M} \in H^1((\wmin,\wmax),\SSS^2)$ satisfying
$\|M+e_3\|_{H^1}, \|\tilde{M}+e_3\|_{H^1} < \delta$,
we have
$$| \Omega(t) - \tilde{\Omega}(t) | \le
C \|M-\tilde{M}\|_{H^{\frac{1}{2}}}, \forall t \in [0,+\infty).$$
Let us prove it only on one of the terms that compose $\Omega$
(the other terms may be treated as well):
$$ \begin{array}{ll}
&
\Big|\int_{\omega_*}^{\omega^*}
\Big( Z' z - \tilde{Z}' \tilde{z}  \Big)
e^{i \omega \int_0^t \zeta} d\omega \Big|
\\
\le &
\Big|\int_{\omega_*}^{\omega^*}
\Big( Z' - \tilde{Z}' \Big) z
e^{i \omega \int_0^t \zeta} d\omega \Big|
\\ & +
\Big|\int_{\omega_*}^{\omega^*}
\tilde{Z}' \Big( z - \tilde{z} \Big)
e^{i \omega \int_0^t \zeta} d\omega \Big|
\\ \le &
\| Z'-\tilde{Z}' \|_{H^{-\frac{1}{2}}}
\| z e^{i \omega \int_0^t \zeta} \|_{H^{\frac{1}{2}}}
\\ & + \| \tilde{Z}' \|_{H^{-\frac{1}{2}}}
\| (z-\tilde{z})e^{i \omega \int_0^t \zeta} \|_{H^{\frac{1}{2}}}
\\ \le &
c_1 K \| M-\tilde{M} \|_{H^{\frac{1}{2}}}. \hfill \Box
\end{array}$$

\noindent \textbf{Proof of Theorem \ref{MainThm}:}
Let $\delta'$ be as in the previous proof.
Let $M^0 \in H^1((\wmin,\wmax),\SSS^2)$ be such that $\| M^0+e_3 \|_{H^1} < \delta'$
and $M \in C^1([0,+\infty),H^1((\wmin,\wmax),\SSS^2))$ be the solution of the closed loop
system such that $M(0)=M^0$.
\\

\noindent \textbf{First step:} Let us prove that $\Omega(t) \rightarrow 0$
when $t \rightarrow + \infty$.} Thanks to the choice of the feedback law,
$M(t)$ is bounded in $H^1$,
uniformly with respect to $t \in [0,+\infty)$.
Computing explicitly $\frac{d\Omega}{dt}(t)$, we see that
$\frac{d\Omega}{dt}(t)$ is bounded in $\mathbb{C}$ uniformly with respect to
$t \in [0,+\infty)-\mathbb{N}T$. Thus, $\Omega$ is uniformly continuous
on $[0,+\infty)$. Since $\Omega \in L^2(0,+\infty)$, it has to satisfy
$\Omega(t) \rightarrow 0$ when $t \rightarrow + \infty$ (Barbalat's lemma).
\\

\noindent \textbf{Second step:} Let us prove that $-e_3$ is the
only possible weak $H^1$ limit.
Let $M_\infty^0$ be a weak $H^1$ limit of the trajectory starting from $M^0$.
There exists a sequence $(t_n)_{n \in \mathbb{N}}$ of $[0,+\infty)$ such that
$t_n \rightarrow +\infty$,
$$\begin{array}{l}
M(t_n) \rightharpoonup M^0_\infty
\text{ weakly in } H^1 \text{ when } n \rightarrow + \infty,
\\
M(t_n) \rightarrow M^0_\infty
\text{ strongly in } H^{1/2} \text{ when } n \rightarrow + \infty.
\end{array}$$
Working as in the previous proof, one may prove that
\begin{equation} \label{H1bound:M(tn)}
\|M(t_n)+e_3\|_{H^1} < \delta, \forall n \in \mathbb{N}.
\end{equation}
There exists $t_\infty \in [0,T)$ such that $t_n \text{ mod } T \rightarrow t_\infty$.
Let $M_\infty(t,\omega)$ be the solution of the closed loop system
associated to the initial condition
$M_\infty(t_\infty)=M_\infty^0$. Let us prove that $\mathcal{L}$ is
constant along this trajectory, by proving that the associated control
$\Omega_\infty$ vanishes. In order to simplify, we assume that $t_\infty=0$
(otherwise, consider an additional shift).
For every $t>0$,
$M(t_n+t) \rightarrow M_\infty(t)$ strongly in $H^{1/2}$ when $n \rightarrow + \infty$,
thanks to the previous proposition.
This allows to pass to the limit in the feedback law:
$\Omega(t_n+t) \rightarrow \Omega_\infty(t)$ when $n \rightarrow + \infty$,
for every $t>0$. Thanks to the first step, we get $\Omega_\infty=0$.

In order to apply  Proposition \ref{Prop:Inv}, we only need to check that
$\|M_\infty^0+e_3 \|_{H^1} < \delta$, which is a consequence of
(\ref{H1bound:M(tn)}). $\hfill \Box$

%-----------------------------------------------------------------------------------------
\subsection{Obstructions to global stabilization}
%-----------------------------------------------------------------------------------------

Now, let us explain why these feedback laws may not provide global stabilization in
$H^1((\omega_*,\omega^*),\mathbb{S}^2)$. The first obstuction is a topological one: the space
$H^1((\omega_*,\omega^*), \mathbb{S}^2)$ cannot be continuously deformed to one point
(because $\mathbb{S}^2$ is not), thus global stabilization in this space is impossible.

Actually, for our explicit feedback laws, it is easy to see that $M^0 \equiv +e_3$
is an invariant solution. It is interesting to know whether it is the only one
(i.e. if one may expect the stabilization of any initial condition $M^0 \neq e_3$).
The answer is no, as emphasized in the following proposition.

\begin{prp} \label{Prop:Invariant-continuum}
For every $\omega_*,\omega^* \in \mathbb{R}$ such that $\omega_*<\omega^*$,
there exists an infinite number of non trivial functions in the LaSalle invariant set.
\end{prp}

The proof is detailed in Appendix \ref{App:Prop-Inv}.
Actually, all the invariant solutions may be computed explicitly.

%OOOOOOOOOOOOOOOOOOOOOOOOOOOOOOOOOOOOOOOOOOOOOOOOOOOOOOOOOOOOOOOOOOOOOOOOOOOOOOOOOOOOOOO
\section{Conclusion} \label{conclusion:sec}
%OOOOOOOOOOOOOOOOOOOOOOOOOOOOOOOOOOOOOOOOOOOOOOOOOOOOOOOOOOOOOOOOOOOOOOOOOOOOOOOOOOOOOOO

We have investigated here the stabilization  of an infinite dimensional system admitting
a continuous spectrum. We have designed  a Lyapunov based feedback.
Closed-loop simulations illustrate  the asymptotic convergence towards the goal steady-state.
We have provided a local and weak convergence  result for the $H^1$ topology. Simulations
indicate that the domain of attraction  is far from being local  and thus we can expect
a large attraction domain for this feedback law. However, the stabilization is not
global because there exists non trivial invariant solutions.

Few problems are still open concerning this problem.
Are the invariant solutions unstable?
Does the local stabilization hold for the strong $H^1$-topology (not only the weak one)?
Is it possible to get semi-global stabilization?
What is the value of convergence rates?
Is it possible to produce arbitrarily fast stabilization?

More generally, this feedback and convergence analysis opens the way to  asymptotic
stabilization of neutrally stable systems of infinite dimension with continuous spectra.
For example, it will be interesting to see if the following  system (1D Maxwell-Lorentz
model for the propagation of an electro-magnetic wave in a non-homogeneous dispersive material)
can also be stabilized to zero:
$$\begin{array}{l}
\dvv{E}{t} + \dvv{P}{t} = \dvv{E}{x},\quad x\in(0,1)
\\
\dvv{P}{t} = p^2(x)(E-P),\quad x\in(0,1)
\\
E(0,t)=u(t),\quad  E(1,t)=v(t)
\end{array}$$
with two controls $u$ and $v$.
When  $p(x)$ is a smooth strictly increasing  positive function,
the above system admits as continuous spectrum  $\pm \imath]p(0),p(1)[$.

%OOOOOOOOOOOOOOOOOOOOOOOOOOOOOOOOOOOOOOOOOOOOOOOOOOOOOOOOOOOOOOOOOOOOOOOOOOOOOOOOOOOOOOOOOOOOOOOOOOOOOOO

%OOOOOOOOOOOOOOOOOOOOOOOOOOOOOOOOOOOOOOOOOOOOOOOOOOOOOOOOOOOOOOOOOOOOOOOOOOOOOOOOOOOOOOOOOOOOOOOOOOOOO
\appendix

%-----------------------------------------------------------------------------------------------------
\section{Proof of Proposition~\ref{Solutions}}  \label{proof:ap}
%-----------------------------------------------------------------------------------------------------

Let $M^0 \in H^1((\wmin, \wmax), \SSS^2)$ and $R>0$ be such that
\begin{equation} \label{def:R}
R> \max\left\{  \|M^0\|_{H^1} , \sqrt{2 \mathcal{L}(0) + \wmax - \wmin } \right\}
\end{equation}
Let $C_1, C_2>0$ be such that
\begin{equation} \label{def:C1}
\begin{array}{l}
\| f e^{-i\omega t} \|_{H^1} \le C_1 \|f\|_{H^1},
\\ \forall f \in H^1(\wmin,\wmax), \forall t \in [0,T],
\end{array}
\end{equation}
\begin{equation} \label{def:C2}
\begin{array}{l}
\| F(t,M_1)-F(t,M_2) \|_{H^1} \le C_2 \|M_1-M_2\|_{H^1},
\\
\forall M_1, M_2 \in B_R[H^1((\wmin,\wmax),\mathbb{R}^3)],
\forall t \in [0,T],
\end{array}
\end{equation}
where $B_R[X]$ denote the closed ball centered at $0$ with radius $R$,
of the space $X$. Let $T^*=T^*(R)>0$ be small enough so that
\begin{equation} \label{def:T*}
\|M^0\|_{H^1} + T^* C_1 K_p (G+2T)R^3 < R
\text{  and  }
T^* C_2 < 1.
\end{equation}
Let us consider the map $\Theta$,
defined on the space
$$E:=B_R[C^0([0,T^*],H^1((\wmin,\wmax),\mathbb{R}^3))]$$
by
$$\Theta(M)(t,\omega):=M^0(\omega)+\int_0^t F(s,M(s,\omega)) ds$$
for every $(t,\omega) \in [0,T^*] \times (\wmin,\wmax)$.

\emph{First step: Let us prove that $\Theta$ takes values in $E$.}
Let $M \in E$. It is clear that $\Theta(M)$ is continuous in time
with values in $H^1((\wmin,\wmax),\mathbb{R}^3)$.
For $t \in [0,T^*]$, we have
$$\| \Theta(M)(t) \|_{H^1}
\le \|M^0\|_{H^1} + \int_0^t \| F(s,M(s)) \|_{H^1} ds.$$
By definition, we have
$$\begin{array}{ll}
&\|F(s,M(s))\|_{H^1}^2
\\
= &  \| \Omega(s) z(s) e^{-\imath \omega \int_{0}^s \varsigma} \|_{H^1}^2
+ \| \Re[ \Omega(s) Z(s) e^{-\imath \omega \int_{0}^s \varsigma}  ] \|_{H^1}^2
\\
\le &
| \Omega(s) |^2 C_1^2 ( \|z(s)\|_{H^1}^2 + \|Z(s)\|_{H^1}^2 )
\\
= &
C_1^2  |\Omega(s)|^2  \|M(s)\|_{H^1}^2.
\end{array}$$
Moreover, the Cauchy-Schwarz inequality gives
$$|\Omega(s)| \le K_p(G+2)\|M(s)\|_{H^1}^2,$$
thus, thanks to (\ref{def:T*}), we have
$$\begin{array}{ll}
& \| \Theta(M)\|_{L^\infty((0,T^*),H^1)}
\\ \le  &
\|M^0\|_{H^1} + T^* C_1  K_p(G+2) R^3 \le R.
\end{array}$$

\emph{Second step: Let is prove that $\Theta$ is a contraction.}
For $M_1,M_2 \in E$ and $t \in [0,T^*]$, using (\ref{def:C2}), we get
$$\begin{array}{ll}
& \| \Theta(M_1)(t) - \Theta(M_2)(t) \|_{H^1}
\\ \le & \int_0^t \| F(s,M_1(s))-F(s,M_2(s)) \|_{H^1}
\\ \le & t C_2 \|M_1 - M_2 \|_{L^\infty((0,T^*),H^1)},
\end{array}$$
thus $\Theta$ is a contraction, thanks to (\ref{def:T*}).

\emph{Third step: Let us prove the existence and uniqueness of
strong solutions, defined on $[0,+\infty)$.}
Thanks to the Banach fixed point theorem,
the map $\Theta$ has a unique fixed point.

We have proved that, for every $R>0$,
there exists $T^*=T^*(R)>0$ such that,
for every $M^0 \in B_R[H^1((\wmin,\wmax),\SSS^2)]$,
there exists a unique weak solution
$$M \in C^0([0,T^*],H^1((\wmin,\wmax),\mathbb{R}^3))$$
in the sense
$$\begin{array}{l}
M(t,\omega)=M^0(\omega)+\int_0^t F(s,M(s,\omega)) ds,
\\ \text{ in } H^1(\wmin,\wmax), \forall t \in [0,T^*].
\end{array}$$
From this equality, we deduce that
$$M \in C^1_{pw}([0,T^*],H^1((\wmin,\wmax),\mathbb{R}^3))$$
and
$$\begin{array}{l}
\frac{dM}{dt}(t,\omega)=F(t,M(t,\omega))
\\ \text{ in } H^1((\wmin,\wmax),\mathbb{R}^3),
\forall t \in [0,T^*]-\mathbb{N}T.
\end{array}$$
Since $H^1 \subset C^0$, we also have
$$\begin{array}{l}
\frac{dM}{dt}(t,\omega)=F(t,M(t,\omega))
\\ \forall t \in [0,T^*]-\mathbb{N}T, \forall \omega \in (\wmin,\wmax).
\end{array}$$
This has 2 consequences:
\begin{itemize}
\item $M(t,.)$ takes values in $\SSS^2$ for every $t \in [0,T^*]$,
indeed, $M^0$ does and the following computation is licit
for every $t \in [0,T^*]-\mathbb{N}T$,
$$\frac{d}{dt} \|M(t,\omega)\|^2
= 2 \langle M(t,\omega) , F(t,M(t,\omega)) \rangle =0,$$
\item the computations (\ref{dL/dt=}), (\ref{dL/dt=bis}) are licit,
thus $\mathcal{L}(t)$ is not increasing.
\end{itemize}
Therefore, we have
$$\begin{array}{ll}
& \| M(T^*) \|_{H^1}^2
\\ = & \int\limits_{\wmin}^{\wmax}
 |Z'(T)|^2 + z'(T)^2 + |Z(T)|^2 + z(T)^2 d \omega
\\  \le &
2 \mathcal{L}(T) + \wmax - \wmin
\\ \le &
2 \mathcal{L}(0) + \wmax - \wmin
\le  R^2
\end{array}$$
thanks to (\ref{def:R}).
Thus, we can apply the previous result with $M^0$ replaced by $M(T^*)$:
it provides a solution on $[0,2T^*]$. Iterating this again, we get a solution
defined for every $ t \in [0,+\infty)$. $\hfill \Box$

%-------------------------------------------------------------------------------------------
\section{Proof of Proposition \ref{Prop:Invariant-continuum}} \label{App:Prop-Inv}
%-------------------------------------------------------------------------------------------

%- - - - - - - - - - - - - - - - - - - - - - - - - - - - - - - - - - - - - - - - - - - - -
\subsubsection{Heuristic}
%- - - - - - - - - - - - - - - - - - - - - - - - - - - - - - - - - - - - - - - - - - - - -

In this section, we perform an heuristic to prove that the LaSalle invariant set
contains non trivial solutions. In order to simplify the  notations, we take $G=1$.
\\

Let $M=(x,y,z)$ be in the LaSalle invariant set.
In view of (\ref{CNS-inv}) and after integrating by parts,  we have
\begin{equation} \label{EDO:inv}
\left\lbrace \begin{array}{l}
xz''-zx''=x \text{ on } (\wmin,\wmax),\\
yz''-zy''=y \text{ on } (\wmin,\wmax),\\
x^2+y^2+z^2=1 \text{ on } (\wmin,\wmax),\\
xz'=x'z \text{ and } yz'=y'z \text{ at } \wmin, \wmax.
\end{array}\right.
\end{equation}

\noindent \textbf{First step:} We proceed to eliminations
in order to get an ordinary differential equation involving only $z$.
Differentiating the third equality of (\ref{EDO:inv}), we get
\begin{equation} \label{Dnorm=0}
xx'+yy'+zz'=0 \text{ on } (\wmin,\wmax),
\end{equation}
\begin{equation} \label{DDnorm=0}
xx''+(x')^2+yy''+(y')^2+zz''+(z')^2=0 \text{ on } (\wmin,\wmax).
\end{equation}
Multiplying this equality by $z$ and using
the 3 first equalities of (\ref{EDO:inv}), we get
\begin{equation} \label{(x')2}
z''+z^2-1+z [ (x')^2 + (y')^2 + (z')^2 ] =0 \text{ on } (\wmin,\wmax).
\end{equation}
Derivating this equality, using again the 2 first equalities of (\ref{EDO:inv}) together with
(\ref{(x')2}) and (\ref{Dnorm=0}) we get
$$z z'''+3z^2 z' + z'(1-z'')=0 \text{ on } (\wmin,\wmax).$$
At the points $\wmin$ and $\wmax$, thanks to (\ref{Dnorm=0}) and
the 4th equality of (\ref{EDO:inv}), we have
$$0=z(xx'+yy'+zz')=z'(x^2+y^2+z^2),$$
thus
$$z'=0 \text{ at } \wmin \text{ and } \wmax.$$
We deduce from the 4th equality of (\ref{EDO:inv}) that
$$zx'=zy'=0 \text{ at } \wmin \text{ and } \wmax.$$
Thus, the relation (\ref{(x')2}) provides
$$z''+z^2-1=0 \text{ at } \wmin \text{ and } \wmax.$$
Therefore, if $M=(x,y,z)$ is in the invariant set,
then $z$ solves the following boundary value problem
\begin{equation} \label{EDO:invz}
\left\lbrace \begin{array}{l}
z z'''+3z^2 z' + z'(1-z'')=0 \text{ on } (\wmin,\wmax),\\
z'=z''+z^2-1=0 \text{ at } \wmin \text{ and } \wmax,
\end{array} \right.
\end{equation}
and $(x,y)$ are solutions of the first order system
\begin{equation} \label{EDO:xy}
\left\lbrace \begin{array}{l}
xx'+yy'=-zz' \text{ on } (\wmin,\wmax),\\
(x')^2 + (y')^2=\frac{1-z^2-z''}{z} - (z')^2 \text{ on } (\wmin,\wmax),\\
zx'=zy'=0 \text{ at } \wmin, \wmax.
\end{array} \right.
\end{equation}

\noindent \textbf{Second step:} Let us solve the equation (\ref{EDO:invz}).
We introduce the function $F:=z''+3z^2$. The first equality of (\ref{EDO:invz})
allows to prove that
$$z F' = (F-1) z \text{ on } (\wmin,\wmax).$$
Thus, there exists $C \in \mathbb{R}$ such that $(F-1)=Cz$, i.e.
\begin{equation} \label{EDOz_ordre2}
z''=-3z^2+Cz+1 \text{ on } (\wmin,\wmax).
\end{equation}
Thanks to this equation, we deduce from the second equality of (\ref{EDO:invz}) that
\begin{equation} \label{z(C-2z)}
z(C-2z) \text{ at } \wmin \text{ and } \wmax.
\end{equation}
Multiplying (\ref{EDOz_ordre2}) by $z'$ and integrating over $(\omega_*,\omega)$, we get
$$\frac{1}{2} (z')^2 = - z^3 + \frac{C}{2} z^2 + z + \text{cst}.$$
The left hand side vanishes at the boundary,
and the right hand side is equal to $z$ at the boundary thanks to (\ref{z(C-2z)}),
thus
$$z(\wmin)=z(\wmax).$$
The conclusion of this second step is the existence of a constant $C \in \mathbb{R}$ such that
$z$ solves one of the following systems
$$(\Sigma_1) \left\lbrace \begin{array}{l}
z''=-3z^2+Cz+1 \text{ on } (\wmin,\wmax),\\
z=z'=0 \text{ at } \wmin, \wmax,
\end{array} \right.$$
$$ \left\lbrace \begin{array}{l}
z''=-3z^2+Cz+1 \text{ on } (\wmin,\wmax),\\
z-C/2=z'=0 \text{ at } \wmin, \wmax.\\
\end{array} \right.$$

\noindent \textbf{Third step:} We prove that $(\Sigma_1)$ has admissible solutions
for arbitrarily small intervals $(\omega_*,\omega^*)$.
Multiplying the first equality of $(\Sigma_1)$ by $z'$ and integrating
over $(\omega_*,\omega)$, we get
\begin{equation} \label{z'2}
(z')^2 = -2z^3+Cz^2+2z \text{ on } (\wmin,\wmax).
\end{equation}
The function $t \mapsto -2t^3+Ct^2+2t$ vanishes at $t=0$,
$t=\beta_C:=(C+\sqrt{C^2+16})/4$ (simple roots) and is positive on
$(0,\beta_C)$ thus one may define
$$G(x):=\int\limits_0^x \frac{dt}{\sqrt{-2t^3+Ct^2+2t}},
\forall x \in [ 0 , \beta_C ],$$
$$\alpha_C:=G\left( \frac{C+\sqrt{C^2+16}}{4} \right).$$
Then $G \in C^\infty(0,\beta_C)$,
$G$ is increasing from $0$ to $\alpha_C$ and it has
an infinite derivative at $x=0$ and $x=\beta_C$.
Thus, $G^{-1}:[0,\alpha_C] \rightarrow [0,\beta_C]$
belongs to $C^2(0,\alpha_C)$, it is increasing from
$0$ to $\beta_C$ and its derivative vanishes at $0$ and $\alpha_C$.
Then the function $z:(0,2\alpha_C) \rightarrow \mathbb{R}$,
symmetric with respect to $\alpha_C$ and such that
\begin{equation} \label{def:z}
z(\omega):=G^{-1}(\omega), \forall \omega \in [0,\alpha_C]
\end{equation}
is a solution of $(\Sigma_1)$ with $\omega_*=0$ and
$\omega^*=2\alpha_C$. Let us emphasize that this solution is
admissible, when $C<0$ because it takes values in $[0,\beta_C]$
which in included in $[0,1)$. Notice that
$$\alpha_C \xrightarrow[C \rightarrow - \infty]{} 0,$$
thus, we have built admissible solutions of $(\Sigma_1)$
for arbitrarily small intervals $(\omega_*,\omega^*)=(0,2\alpha_C)$.
\\

\noindent \textbf{Fourth step:} We prove that, for any solution of $(\Sigma_1)$
with $C<0$, the system (\ref{EDO:xy}) has solutions.
First, notice that, thanks to $(\Sigma_1)$ and (\ref{z'2}), we have
$$\frac{1-z^2-z''}{z} - (z')^2 = 2z^3-C z^2 -C.$$
Eliminating $x'$ in the two first equalities of (\ref{EDO:xy}), we get
$$a (y')^2 + by'+c=0$$
where
$$\begin{array}{l}
a:=x^2+y^2=1-z^2,\\
b:=2zz'y,\\
c:=(zz')^2-x^2(2z^3-Cz^2-C).
\end{array}$$
Thanks to (\ref{z'2}) one may prove that the discriminant is
$$\Delta:=b^2-4ac=-4Cx^2.$$
Therefore,
$$y'=-\frac{zz'}{1-z^2} y \pm \frac{\sqrt{|C|}x}{1-z^2}.$$
By symmetry, we also have
$$x'=-\frac{zz'}{1-z^2} x \pm \frac{\sqrt{|C|}y}{1-z^2}.$$
In order to ensure $xx'+yy'+zz'=0$, the signs $\pm$ need to be opposite.
Therefore $(x,y)$ are solutions, for instance, of the following linear system
\begin{equation} \label{EDOLxy}
\left\lbrace \begin{array}{l}
x'=-\frac{zz'}{1-z^2} x + \frac{\sqrt{|C|}}{1-z^2} y,\\
y'=-\frac{zz'}{1-z^2} y - \frac{\sqrt{|C|}}{1-z^2} x,\\
x(0)=x_0,\\
y(0)=y_0,
\end{array}\right.
\end{equation}
where $x_0$ and $y_0$ are real numbers such that
$x_0^2+y_0^2=1$.
\\

\noindent \textbf{Conclusion:} For every $C<0$, we have built candidates
of invariant solutions of the closed loop system
associated to the interval $(\wmin,\wmax)=(0,2\alpha_C)$,
which is arbitrarily small when $C \rightarrow - \infty$.
In order to conclude, one just needs to check that
this candidate indeed solves (\ref{EDO:inv}),
which will be done rigourously in the next subsection.

\begin{rmq}
Let us emphasize that the same phenomena happens if we put different gains in the
Lyapunov functions: for every $G_1, G_2 \geqslant 0$, the feedback laws associated to the
control Lyapunov function
$$\mathcal{L}(t):=\int_{\wmin}^{\wmax} \Big[
|Z'|^2 + (z')^2 + G_1 z + G_2 |Z|^2
\Big] d\omega$$
generate a non trivial LaSalle invariant set.
\end{rmq}

%- - - - - - - - - - - - - - - - - - - - - - - - - - - - - - - - - - - - - - - - - - - - -
\subsubsection{Rigorous proof}
%- - - - - - - - - - - - - - - - - - - - - - - - - - - - - - - - - - - - - - - - - - - - -

\begin{lm}
There exists a continuous function
$$\begin{array}{ccl}
(-\infty,0) & \rightarrow & (0,+\infty) \\
     C      & \mapsto     & L_C
\end{array}$$
such that,
\begin{itemize}
\item for every $C \in (-\infty,0)$,
there exists a function $z \in C^3([0,L_C],[0,1))$ such that $z>0$ on $(0,L_C)$ and
$$\left\lbrace \begin{array}{l}
z''=-3z^2+Cz+1 \text{ on } (0,L_C),\\
z=z'=0 \text{ at } 0 \text{ and } L_C,
\end{array}\right.$$
\item $L_C \rightarrow 0$ when $C \rightarrow - \infty$.
\end{itemize}
\end{lm}

\noindent \textbf{Proof:} Consider $L_C:=2\alpha_C$ and
$z:[0,L_C] \rightarrow \mathbb{R}$, symetric with respect to $\alpha_C$
and defined by (\ref{def:z}). $\hfill \Box$

\noindent \textbf{Proof of Proposition \ref{Prop:Invariant-continuum}:}
Let $-\infty < \omega_* < \omega^* < +\infty$.
Let $(N,C) \in \mathbb{N}^* \times \mathbb{R}^*_+$ be such that $\omega^*-\omega_*=NL_C$
(the existence is ensured by the intermediate values theorem and
there exists an infinite number of such couples).
Let $z:[0,L_C] \rightarrow [0,1)$ be as in the previous Lemma.
Let $x_0, y_0 \in \mathbb{R}$ be such that $x_0^2+y_0^2=1$, and
$(x,y)$ be the solution of (\ref{EDOLxy}): this solution is well defined
on the whole interval $[0,L_C]$ because the system is linear and
its coefficients are continuous. Now, let us check that (\ref{EDO:inv}) holds
with $(\omega_*,\omega^*)$ replaced by $(0,L_C)$.

\noindent \textbf{First step:} We check that the 3rd equality of (\ref{EDO:inv}) holds.
Thanks to (\ref{EDOLxy}), the quantity $N:=x^2+y^2+z^2$ solves
$$\frac{d}{dt} ( 1-N )=\frac{zz'}{1-z^2} (1-N) \text{ on } (0,L_C)$$
and $(1-N)(0)=0$, thus $N \equiv 1$ on $(0,L_C)$.
\\

\noindent \textbf{Second step:} We check that the fourth equality of (\ref{EDO:inv}) holds.
Since $z=z'=0$ at $0$ and $L_C$, we also have
$xz'=x'z$ and $yz'=y'z$ at $0$ and $L_C$.
\\

\noindent \textbf{Third step:} We check that the 2 first equations of (\ref{EDO:inv}) hold.
The computations are similar to the ones of the Heuristic
but now, we know that the functions considered are smooth (they are explicit),
so these computations are licit.
Starting from (\ref{EDOLxy}), using the result of the first step, the equality
(\ref{z'2}), and the 2nd order equation solved by $z$, we get
\begin{equation} \label{*}
(x')^2 + (y')^2 +(z')^2 = 2z-C=\frac{1-z^2-z''}{z}.
\end{equation}
Thus, we have
\begin{equation} \label{*bis}
z''+z^2-1+z [ (x')^2 + (y')^2 +(z')^2 ]=0.
\end{equation}
The result of the first step justifies
$$(x^2+y^2+z^2) z''+ [ (x')^2 + (y')^2 +(z')^2 ] z - (x^2+y^2)=0,$$
that may be written
\begin{equation} \label{1}
x(xz''-x)+y(yz''-y)+z^2z''+ [ (x')^2 + (y')^2 +(z')^2 ] z=0.
\end{equation}
Now, differentiating 2 times the identity $x^2+y^2+z^2=1$ and
multiplying the resulting equality by $z$, we get
\begin{equation} \label{2}
xx''z+yy''z+z^2z''+ [ (x')^2 + (y')^2 +(z')^2 ] z=0
\end{equation}
Thus, (\ref{1})-(\ref{2}) gives
\begin{equation} \label{3}
x(xz''-x-x''z)+y(yz''-y-y''z)=0.
\end{equation}
Derivating (\ref{*bis}) and multiplying the resulting equality by $z$, we get
\begin{equation} \label{4}
zz'''+2z^2z'+z'(1-z^2-z'')+2z^2(x'x''+y'y''+z'z'')=0.
\end{equation}
We deduce from the second order equation solved by $z$ that
\begin{equation} \label{4bis}
z z'''+3z^2z'+z'(1-z'')=0.
\end{equation}
Indeed, we have
$$zz'''=z(-6zz'+Cz'),$$
$$3z^2z'+z'(1-z'')=z'[3z^2+1+3z^2-Cz-1]=zz'(6z-C).$$
Thanks to the identity $xx'+yy'+zz'=0$,
the equation (\ref{4bis}) may be written
$$zz'''+2z^2z'+z'(1-z^2-z'')+2z[ z''(xx'+yy'+zz') + zz' ]=0$$
or
\begin{equation} \label{5}
\begin{array}{l}
zz'''+2z^2z'+z'(1-z^2-z'')
\\
+2z[x'(z''x-x)+y'(z''y-y)+z''z'z]=0
\end{array}
\end{equation}
Finally, (\ref{5})-(\ref{4}) gives
\begin{equation} \label{6}
x'(xz''-x-x''z)+y'(yz''-y-y''z)=0.
\end{equation}
Now, (\ref{3}) and (\ref{6}) give the conclusion, because
$$\text{det} \left(
\begin{array}{cc}
x & y \\ x' & y'
\end{array}
\right)
=\frac{\sqrt{|C|}}{1-z^2}(x^2+y^2)=\sqrt{|C|} \neq 0.$$

\noindent \textbf{Conclusion:} We build an invariant solution
$\tilde{M}$ on $(\omega_*,\omega^*)$.
This solution is $L_C$-periodic and satisfies
$\tilde{M}(\omega):=(x,y,z)(\omega-\omega_*)$,
$\forall \omega \in [\omega_*,\omega_*+L_C]$. $\hfill \Box$

\end{document}